\documentclass{article}

\usepackage{amsmath}  
\usepackage{amssymb}  
\usepackage{amsfonts} 
\usepackage{epsfig}
\usepackage{hyperref}
\usepackage{url}
\usepackage{epigraph}

\begin{document}

\title{High-Degree Polynomial Approximations for Solving Linear Integral, Integro-Differential, and Ordinary Differential Equations}

\author{V. V. Kryzhniy }
\date{May 11, 2024}
\maketitle
\begin{abstract}
This paper presents a universal numerical scheme tailored for tackling linear integral, integro-differential, and both initial and boundary value problems of ordinary differential equations. The numerical scheme is readily adapted for resolving ill-posed problems.  Central to our approach is high-degree piecewise-polynomial approximation to the exact solution.  We illustrate the accuracy and stability of our numerical solutions in the presence of noise  through illustrative examples.
Additionally, we demonstrate that proposed regularization being applied  to high-degree interpolation, effectively eliminates Runge's phenomenon.

\end{abstract}
\setlength\epigraphwidth{.5\textwidth}
\setlength\epigraphrule{0pt}

\textbf{Keywords:} Integral equations, numerical methods, inverse problems, ill-posed problems, integro-differential equations, initial and boundary value problems, ordinary differential equations,  numerical interpolation\\
\textbf{MSC codes:}  65D05, 65L05, 65L60, 65R30, 65R32

\section{Introduction}

In this paper, we will present a universal numerical scheme for solving linear equations that include differentiation and integration. For definiteness, let's consider the following integro-differential equation.\footnote{Higher-order derivatives and more integral terms can also be included in the consideration.}

\begin{eqnarray}
\label{4}
g_2(x) y''(x)  + g_1(x) y'(x) + g_0(x) y(x)  + \\ \nonumber h_1(x) \int_a^b {K_1(x, t) y(t)\mathrm{d}t}  + h_2(x) \int_a^x{K_2(x, t) y(t) \mathrm{d}t} = \phi(x), 
\end{eqnarray} 
where $y,  y',  y'' $ are the sought function and its derivatives;

$g, h $ are  functional multipliers;
kernels of integral terms, $K_1(x, t)$ and $ K_2(x, t) $, are integrable functions;  and the right-hand side $\phi(x)$ is given on a set of points  $x = x_0, x_1, \dots, x_{m-1}$.

The numerical scheme is  constructed by  approximating  the sought function  $y(x)$  on  the closed interval $[\alpha, \beta]$ by a piecewise-polynomial with a flexible number of subintervals and degrees of partial polynomials. 

We acknowledge that only specific cases of Equation (\ref{4}) have been extensively studied theoretically. For the purposes of this paper, we assume the existence and uniqueness of a solution to the equation under consideration and that the solution can be approximated by a piecewise-polynomial.

The paper is organized as follows. In the next section, we discuss the piecewise-polynomial approximation of a dataset by solving a minimization problem.

Then, we consider the discretization of Eq. (\ref{4}) with initial/boundary value constraints. In the fourth section, we discuss the regularization of the discretized equation.

Finally, we illustrate the numerical scheme by considering examples of solving various particular cases of Eq. (\ref{4}).

 \section{High-Degree Polynomial Approximation \label{S2} }
 Polynomials are highly convenient functions for discretizing equation (\ref{4}) because they can  be easily evaluated, summed, differentiated and integrated.
 
Moreover,  Weierstrass's approximation theorem guarantees that  any continuous function $f(x)$ on the closed interval  can be approximated by a polynomial    $p_n(x)$, and  $\lim_{n \to \infty}  p_n(x) = f(x)$. Hence, theoretically, high degree polynomial approximations are beneficial.
\\

However, in practice,  high degree polynomials are often avoided  due to ill-conditioning of high degree interpolation  \cite{ Kahaner, Press}.

\begin{figure}[ht]
\begin{center}
\includegraphics[height=9.0 cm]{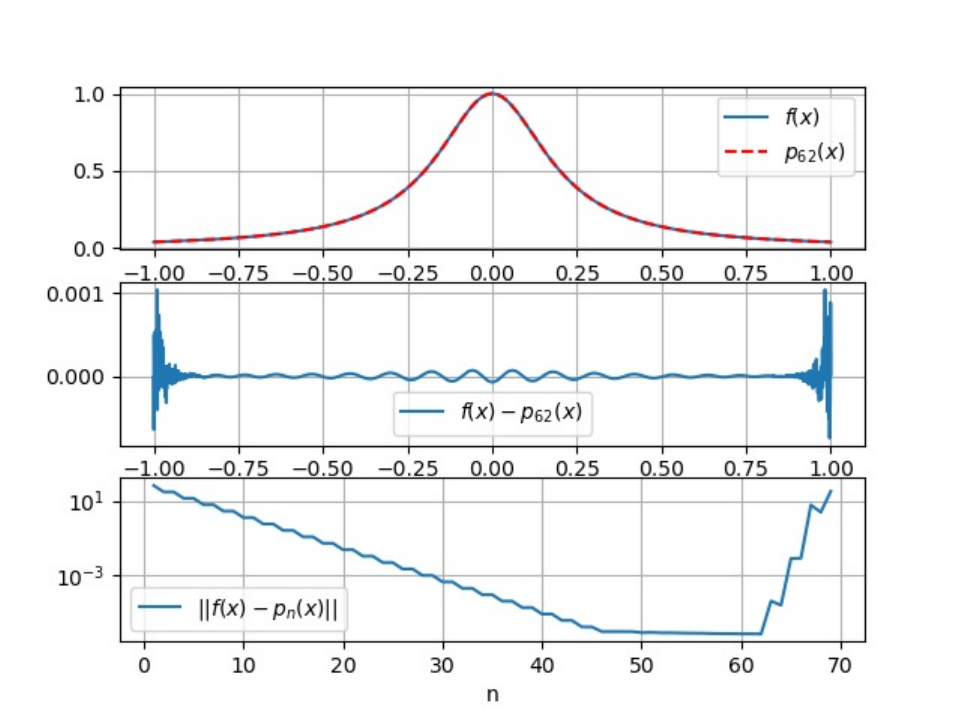}
\end{center}
\caption{Polynomial approximation of Runge's function}
\label{Fig0}
\end{figure}

In modern times, with vast computational power at our disposal, we are not constrained to using interpolating polynomials. 
Assuming that the number of data points significantly exceeds the number of polynomial coefficients, we can determine the coefficients of the polynomial approximation  $ y(x) = \sum_{i = 0}^{n-1} c_i x^i \approx  \phi(x)$  by minimizing the distance between the data and its approximation:
\begin{equation}
\label{7}
 \min_c || A \hat {c} - \hat{\phi} ||^2,
\end{equation} 
where $A$ is a $m \times n$ matrix used for computing polynomial at points $x = x_0, x_1, \cdots, x_{m-1}$.

Then, the best-fit coefficients $\hat{c}$ of a polynomial of a fixed degree are computed as the least squares solution of Eq. (\ref{7}):
\begin{equation}
\label{8}
 \hat{c} = (A^T A)^{-1} A^T \hat{\phi}, 
\end{equation} 
where $A^T$  is a transpose of matrix $A$.

Let us explore polynomial approximations  of Runge's function  $f(x) = (1 + 25x^2)^{-1}$  using polynomials of  increasing degree.  
The polynomial coefficients were computed using  $ m = 800$ uniformly distributed double precision points  on the interval 
$[-1, 1]$.

The results are presented in Figure \ref{Fig0}. The graph in the upper frame displays the best polynomial fit obtained in the experiment, while the graph in the middle frame shows the achieved accuracy. The graph in the lower frame illustrates that the fitting accuracy increases with the polynomial degree up to a certain limit, denoted as $n_{max}$, which is constrained by the fixed-point computer arithmetic.

Therefore, the best polynomial fit of the given data can be found by enumerating all possible polynomial degrees:
 \begin{equation}
\label{9}
\min_{n < n_{max}} \min_{c}|| A \hat {c} - \hat{\phi} ||^2.
\end{equation} 

It is evident, however, that the upper bound on the polynomial degree,  $n_{max}$, limits our ability to find a good polynomial approximation for any continuous function on a finite closed interval $[\alpha, \beta]$. 

Then, we can divide the interval $[ \alpha, \beta]$  into $k $  subintervals and find the best attainable  piecewise-polynomial approximation 
 \footnote{ For simplifying notations, we will denote a piecewise-polynomial approximation as $pp(x)$.} 

\begin{equation}
\label{pp}
y(x) \approx \begin{cases} 
\sum_{i=0}^{n_1} c_i x^i  & \text{for } x \in [\alpha, x_1], \\
\sum_{i=0}^{n_2} c_i x^i  & \text{for } x \in (x_1, x_2], \\
\vdots \\
\sum_{i=0}^{n_{k}} c_i x^i  & \text{for } x \in (x_{k-1}, \beta]
\end{cases}
\end{equation}
by solving the following minimization problem:
 
 \begin{equation}
\label{10}
\min_{k} \min_{n_1 < n_{max}} \dots  \min_{n_k < n_{max}} \min_{c}|| A \hat {c} - \hat{\phi} ||^2, 
\end{equation} 
where $A$ becomes a block-diagonal matrix, and $\hat{c}$ becomes a vertically stacked coefficients of partial polynomials.

\begin{figure}[ht]
\begin{center}
\includegraphics[height=9.0 cm]{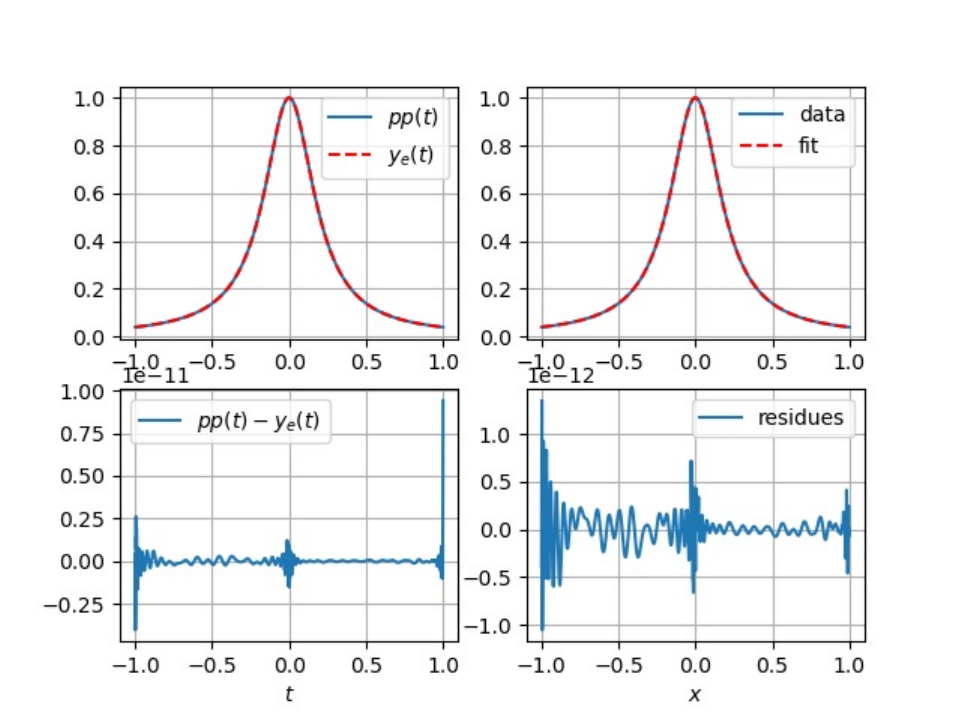}
\end{center}
\caption{Piecewise-polynomial approximation of Runge's function.}
\label{Fig1}
\end{figure}

With two or more subintervals  the minimization  (\ref{10}) is essentially global minimization.
Nonetheless, it yields a highly accurate piecewise-polynomial approximation when the number of data points in each subinterval significantly exceeds the number of polynomial coefficients.

The graph in Figure \ref{Fig1} illustrates the piecewise-polynomial approximation of Runge's function,  as previously considered.

In this and subsequent figures, the graphs in the right-hand column depict the data and fitting residuals, while those in the left-hand column display the exact solution and its piecewise-polynomial approximation, computed at a rate of five points per data point.

 It's worth mentioning that we can introduce constraints to guarantee a desired degree of smoothness of the piecewise-polynomial at joints. In all examples, continuity of the piecewise-polynomial, its derivative, and the second derivative have been enforced.

\section{Discretization and Transformation into a Minimization Problem   \label{S3}}

Assuming that the exact solution of equation (\ref{4}) exists, is unique, and continuous, we aim to find a piecewise-polynomial approximation of the exact solution 
 $y(x)  \approx pp(x)$ on a closed interval $x \in [\alpha, \beta]$.

In cases where equation (\ref{4}) contains derivatives, the unique solution exists if additional boundary or initial conditions are provided:
\begin{eqnarray}
\label{11}
y(\alpha) = y_0,    y'(\alpha) = y_1  \text{- initial  conditions.  } \\
\label{12}
 y(\alpha) = y_{\alpha},    y(\beta) = y_{\beta} \text{ -  boundary conditions.   } 
\end{eqnarray}

 Let  $\hat{\phi} ,   \hat{g}_0 , \hat{ g}_1,  \hat{g}_2 , \hat{h}_1, \text{ and  }  \hat{h}_2  $ denote vectors of length  $ m $ obtained by computing the right-hand side and functional multipliers  of Eq. (\ref{4})  $ \phi(x), g_0(x), g_1(x), g_2(x), h_1(x), h_2(x) $ at points $x _i \in  [\alpha, \beta],  i = 0, 1, \dots, m-1$ ;  $ m  \gg n$ .

Next,  we represent the evaluation of  all  terms in Eq. (\ref{4})  at points  $x _i \in  [\alpha, \beta],  i = 0, 1, \dots, m-1$  in matrix form.

 The evaluation of $y(x),  y'(x), y''(x)$ is given by the following matrix expressions:
  \footnote{For simplicity, we present these matrices only for polynomial approximations over the entire interval.}

 \begin{equation}
\label{13}
 Y \hat{c} = \begin {bmatrix}
1 & x_0^1 & x_0^2 & \cdots & x_0^{n-1} \\
1 & x_1^1 & x_1^2 & \cdots & x_1^{n-1} \\
\vdots & \vdots & \vdots &  \ddots & \vdots \\
1 & x_{m-1}^1 & x_{m-1}^2 & \cdots & x_{m-1}^{n-1} 
\end{bmatrix}
 \begin {bmatrix}
 c_0 \\
 c_1 \\
 \vdots \\
 c_{n-1}
 \end{bmatrix}
\end{equation}

\begin{equation}
\label{14}
Y' \hat{c} =  \begin {bmatrix}
0 & 1 & 2 x_0^1  & \cdots & (n-1)x_0^{n-2} \\
0 & 1 &2 x_1^1 & \cdots & (n-1)x_1^{n-2} \\
\vdots & \vdots & \vdots &  \ddots & \vdots \\
0 & 1 & 2 x_{m-1}^1  & \cdots & (n-1)x_{m-1}^{n-2} 
\end{bmatrix}
 \begin {bmatrix}
 c_0 \\
 c_1 \\
 \vdots \\
 c_{n-1}
 \end{bmatrix}
\end{equation}

\begin{equation}
\label{15}
Y'' \hat{c} = \begin {bmatrix}
0 & 0 & 2 & 6 x_0^1  & \cdots & (n-1)(n-2)x_0^{n-3} \\
0 & 0 &2  & 6 x_1^1 & \cdots & (n-1)(n-2)x_1^{n-3} \\
\vdots & \vdots & \vdots &\vdots & \ddots & \vdots \\
0 & 0 & 2  & 6 x_{m-1}^1  & \cdots & (n-1)(n-2)x_{m-1}^{n-3} 
\end{bmatrix}
 \begin {bmatrix}
 c_0 \\
 c_1 \\
 \vdots \\
 c_{n-1}
 \end{bmatrix}
\end{equation}

The Fredholm and Volterra integral terms will become:
\begin{equation}
\label{16}
M \hat{c }= \begin {bmatrix}
\mu_{0,0}& \mu_{0,1} & \mu_{0,2}  & \cdots  &\mu_{0,(n-1)} \\
\mu_{1,0}& \mu_{1,1} & m_{1,2}   &\cdots  & \mu_{1,(n-1)} \\
\vdots & \vdots & \vdots  &\ddots & \vdots \\
\mu_{(m-1),0}& \mu_{m,1} & \mu_{(m-1),2}  & \cdots  & \mu_{(m-1),(n-1)} \\
\end{bmatrix}
 \begin {bmatrix}
 c_0 \\
 c_1 \\
 \vdots \\
 c_{n-1}
 \end{bmatrix},
\end{equation}
 where $\mu_{i,j} = \int_a^b {K(x_i t) t^j} \mathrm{d} t $,

\begin{equation}
\label{17}
V \hat{c }= \begin {bmatrix}
v_{0,0}& v_{0,1} & v_{0,2}  & \cdots  & v_{0,(n-1)} \\
v_{1,0}& v_{1,1} & v_{1,2}  & \cdots & v_{1,(n-1)} \\
\vdots & \vdots & \vdots & \ddots & \vdots \\
v_{(m-1),0}& v_{m,1} & v_{(m-1),2}  & \cdots  & v_{(m-1),(n-1)} \\
\end{bmatrix}
 \begin {bmatrix}
 c_0 \\
 c_1 \\
 \vdots \\
 c_{n-1}
 \end{bmatrix},
\end{equation}
 where $v_{i,j} = \int_a^{x_i} {K(x_i t) t^j} \mathrm{d} t$.
 
 By denoting element-wise multiplication of each column of a matrix by a column vector as  a dot, we  obtain Eq. (\ref{4}) in matrix form:
\begin{equation}
\label{18}
(g_2\cdot Y'') \hat{c} + (g_1\cdot Y' )\hat{c} +( g_0\cdot Y)\hat{ c} + (h_1\cdot M)\hat{ c} + (h_2 \cdot V )\hat{c} =  A \hat{c} = \hat{\phi}
\end{equation}

When  the initial conditions (\ref{11})  are provided, they will hold  true if the polynomial's coefficients are such that:
\begin{equation}
\label{19}
\left\{
\begin{array}{c}
1  c_0   +  \alpha c_1  + \alpha^2 c_2  + \cdots  +  \alpha^{n-1}  c_{n-1} = y_0 \\
0 c_0    +     1 c_1      +  2 \alpha  c_2  + \cdots  +  (n-1) \alpha^{n-2}  c_{n-1} = y_1
\end{array}
\right.
\end{equation} 

Similarly, for boundary conditions the constraints are:

\begin{equation}
\label{20}
\left\{
\begin{array}{c}
1  c_0   + \alpha c_1  + \alpha^2 c_2  + \cdots  +  \alpha^{n-1} c_{n-1} = y_{\alpha} \\
1  c_0  +  \beta  c_1      +  \beta^2 c_2  + \cdots  +   \beta^{n-1}  c_{n-1} = y_{\beta}
\end{array}
\right.
\end{equation} 

It follows from the constraints above that there are $n-2$ independent polynomial's coefficients. 
 Denoting these independent coefficients as $\tilde{ c}$, we can express all  polynomial coefficients in terms of these independent ones as:
\begin{equation}
\label{21}
\hat{c} = T  \tilde{ c} + {s},
\end {equation}
where $T$ is a matrix, and $s$ is a vector.

Let's mention that in the case of piecewise-polynomial approximation, the equality constraints (\ref{19}) or (\ref{20}) will include smoothness constraints at joints of polynomial parts. 

By substituting (\ref{21}) into (\ref{18}) and denoting $\tilde{A} =  A  T,  \tilde{ \phi} = \hat{\phi}  - A s  $ we can represent  equation (\ref{4}) with initial (\ref{11}) or boundary conditions  (\ref{12})  in matrix form:

\begin{equation}
\label{22}
\tilde{A}  \tilde{c} =  \tilde{\phi}
\end {equation}

Then, the independent coefficients of the  best-fit polynomial are computed by solving the minimization problem:
\begin{equation}
\label{23}
 \min_{\tilde{c}}|| \tilde{A}  \tilde{c} - \tilde{\phi} ||^2
\end{equation} 
as described in the previous section.

Finally, the sought function is obtained by computing coefficients $\hat{c}$ using  (\ref{21}) and evaluating the polynomial at any set of points 
$t_i \in [\alpha, \beta]$.

\section{Computing Regularized Polynomial Coefficients \label{S4}}
 
 Ill-posed integral equations of the first kind, which require regularization, are significant particular cases of Equation (\ref{4}).
  Let's begin with the quadrature discretization method commonly used  for regularizing the Fredholm integral equation of the first kind:
  \begin{equation}
\label{1}
\int_{\alpha}^{\beta} {K(x, t) y(t) \mathrm{d}t }= \phi(x) . 
\end{equation}
  
By computing integral in (\ref{1})  with the help of  a $n$-point quadrature rule, the integral equation is transformed into a matrix equation:
 \begin{equation}
\label{24}
A \hat{y} = \hat{\phi}, 
\end {equation}
where $ A $ is a $ m \times n$ matrix, and $\hat{y}, \hat{\phi} $ are vectors of length $n, m$, respectively.

Then, according to Tikhonov's theory \cite{Tikhonov}, a regularized solution is obtained by solving  the following minimization problem:

\begin{equation}
\label{25}
 \min \left\{ || A  \hat{y} - \hat{\phi} ||^2 + \lambda || L \hat{y} ||^2 \right\},
\end {equation}
where $ \lambda$ is the regularization parameter, and $L $  is a stabilizing matrix.

For given $\lambda > 0$ and $L$ the least squares solution of problem (\ref{25}) will be \cite{Hansen2}:
\begin{equation}
\label{26}
y_{\lambda} = (A^TA + \lambda L^TL)^{-1}A^T \hat{\phi}
\end {equation}

The stabilizing matrix $L$ is typically either the identity matrix or a matrix that approximates the second derivative of the solution \cite{Hansen2}.
The optimal value  of the regularization parameter $\lambda$ is computed using a special criterion  \cite{ Hansen2, Tikhonov}.
\\

Now, let's  discretize the Eq. (\ref{1}) by approximating the solution with a polynomial  $ y(x)  \approx \sum_{i=0}^{n-1}  c_i x^i$ .
Following the notations from the previous section,  the corresponding matrix equation is:

 \begin{equation}
\label{27}
M\hat{c} = \hat{\phi}.
\end {equation}

Similarly to standard regularization, we can formulate  a  minimization problem:
 
 \begin{equation}
\label{28}
 \min_{\hat{c}} \left\{ || M  \hat{c} - \hat{\phi} ||^2 + \lambda || L \hat{c} ||^2 \right\}, 
\end {equation}

Then, the regularized polynomial coefficients $\hat{c}$  are  computed in a manner identical to  (\ref{26}):

\begin{equation}
\label{29}
\hat{c}_{\lambda} = (M^TM + \lambda L^TL)^{-1}M^T \hat{\phi}
\end {equation}

Regularization of Eq. (\ref{27}) is somewhat nonstandard because we lack a priori information or expected properties for the coefficients $c_i$.

It is important to stress that the matrix equation (\ref{27}), obtained for a fixed number of components, is traditionally solved without regularization using methods such as collocation, least squares, or the Galerkin method \cite{Polyanin}.
\\

We can regularize equation (\ref{27})  by selecting  $L=Y$   (\ref{13})  to  limit the norm of the approximate solution or  $L=Y''$ (\ref{15}) to limit its second derivative.

Considering that only the product $L^TL$ is present in formula (\ref{29}), we can compute matrices $Y$ and $Y''$ at any desirable set of points 
$x_i \in [\alpha,\beta],  \text{ where } i = 0, 1, \dots (k-1)$. In particular, we can take $k \gg m$.

It's evident that the constructed regularizing operator depends  on the parameters of the piecewise-polynomial approximation and the regularization parameter. 
According to the conjecture proposed in \cite{kr2}, optimal values of all parameters can be obtained using  a criterion to find  the regularization parameter. 
Thus, in the case of the generalized cross validation criterion \cite{Golub, Hansen2},  the minimization problem for finding coefficients of the best fit polynomial (\ref{9})  will be replaced with: 
\begin{equation}
\label{30}
\min_{n < n_{max}} \min_{c, \lambda}{\frac{||M \hat{c}_{\lambda} -  \hat{\phi}||^2}{\mathrm{trace}(I - MM^{\sharp})^2}},
\end{equation} 
where  $ M^{\sharp} = (M^TM + \lambda L^TL)^{-1}M^T $,  and $I$ is the identity matrix.

\section{Solution of Well-Posed Problems }
In this section we illustrate  the proposed numerical schema by solving ordinary differential equations with initial or boundary conditions, integro-differential equations, and inhomogeneous  integral equations on the second kind. In all examples we assume that the number of data point exceeds the total number of polynomial coefficients $ m \gg n$.  

\subsection{Linear ordinary differential equations \label{S52}}
As shown in Section \ref{S3}, there are no differences in solving differential equations with initial or boundary conditions.
Let's solve the boundary value problem for the Bessel differential equation of order $\nu = 0 $ on the interval $[0, 100]$.

\begin{eqnarray}
\label{31}
x^2  y''(x)  + x y'(x) + (x^2 - \nu^2 )y(x)  = 0,   \\
\nonumber  \nu = 0 \\
 \nonumber y(0) = J_0(0)\\
 \nonumber   y(100) = J_0(100) 
\end{eqnarray} 

The  Bessel differential equation (\ref{31}) has a regular singular point at $ x = 0$. The exact   solution of problem(\ref{31}) is the Bessel function $y(x) = J_0(x)$.

According to discussion in Section \ref{S3}, the task (\ref{31}) is transformed into a minimization problem (\ref{23}). Then, for a fixed polynomial degree, the least square solution of the latter is computed using formula (\ref{8}).  It can be observed from equation (\ref{8}) that the solution of the homogeneous equation is zero. 

In  cases where it is known that the task under consideration has a unique solution, we can use a simple change of variables to rewrite the identical non-homogeneous problem.

Let's rewrite task (\ref{31}) for the function $z(x) = y(x) + 1$. 
\begin{eqnarray}
\label{32}
x^2  z''(x)  + x z'(x) + x^2 z(x)  = x^2,   \\
 \nonumber z(0) = J_0(0) + 1\\
 \nonumber   z(100) = J_0(100) + 1
\end{eqnarray} 

\begin{figure}[ht]
\begin{center}
\includegraphics[height=9.0 cm]{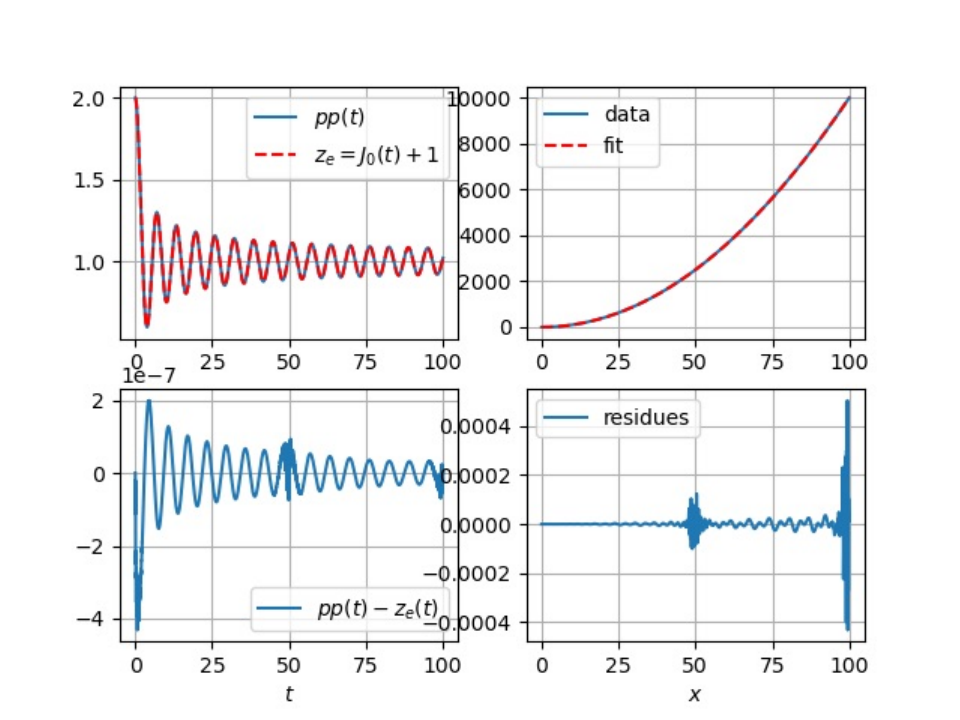}
\end{center}
\caption{ Solution of boundary value problem for Bessel ODE  in form (\ref{32})}
\label{Fig5}
\end{figure}

For numerical calculations, the right-hand side of Eq. (\ref{32}) has been computed at 
$ m=1000$  double precision points uniformly distributed on the interval $[0,100]$ . The piecewise-polynomial approximation 
$pp(t)$  of the solution is shown in the graph in Figure \ref{Fig5}. 

Next, let's consider the solution of an initial value problem:
\begin{eqnarray}
\label{33}
y''(x)  + y'(x) + y(x)  = \cos x  - \sin x   + \delta \\ 
 \nonumber y(0) =  1\\
 \nonumber   y'(0) =   1,
\end{eqnarray} 
where $\delta$ is a normally distributed noise with $ \sigma = 0.01$. 

\begin{figure}[ht]
\begin{center}
\includegraphics[height=9.0 cm]{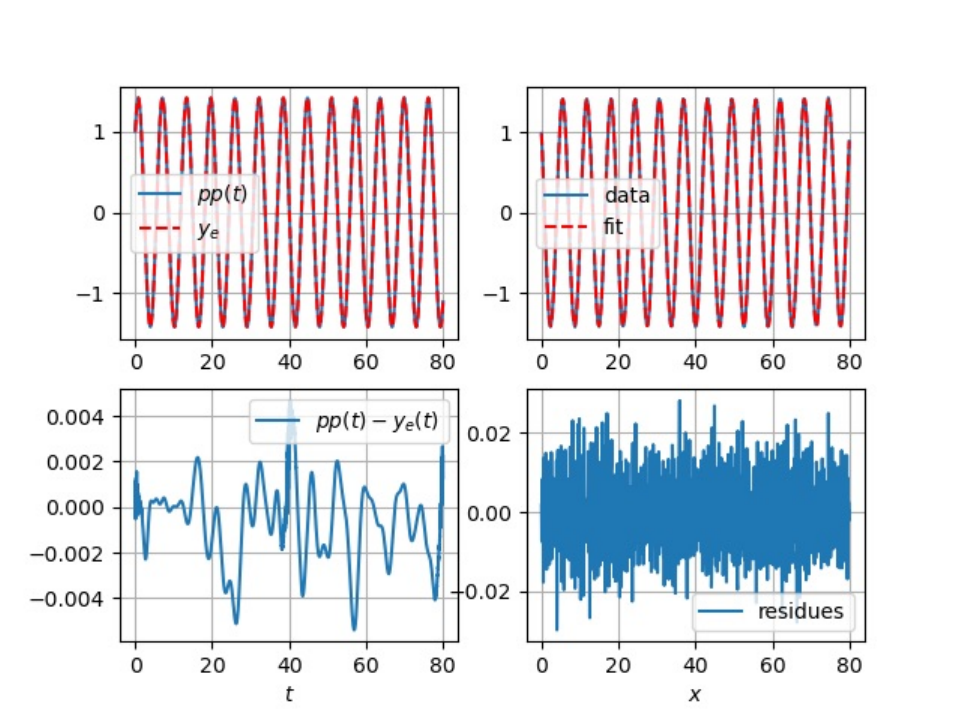}
\end{center}
\caption{ Solution of initial value problem (\ref{33}). }
\label{Fig6}
\end{figure}

The exact solution of problem (\ref{33}) is $y_e(x) = \sin x + \cos x$.  As can be seen from the graph in Fig. \ref{Fig6}, the solution computed using $m=800$ uniformly distributed noisy points is even more accurate than the right hand side of Eq. (\ref{33}).

\subsection{Linear integro-differential equations \label{S53}}

Since there is a vast variety of integro-differential equations (IDEs), the analytical and numerical methods for solving them vary considerably.  
The presented numerical schema allows for solving IDEs using a universal approach.

Let's consider a Volterra-type integro-differential equation.
\begin{eqnarray}
\label{34}
y'(x) - 2 \int_0^x \sin(x - t) y(t) \mathrm{d}t = \cos{x}  + \delta  \\
\nonumber y(0) =  1
 \end{eqnarray} 
 
The exact solution $y(x)$ can be found using Laplace transform method.
 
\begin{equation}
\nonumber  y(x) = \frac{3 \mathrm{e}^x }{4}  + \frac{\mathrm{e}^{-x/2} (2\omega  \cos{\omega x} + 3 \sin{\omega x})}{8 \omega},    \omega = \sqrt{7}/2
 \end{equation} 
 
\begin{figure}[ht]
\begin{center}
\includegraphics[height=9.0 cm]{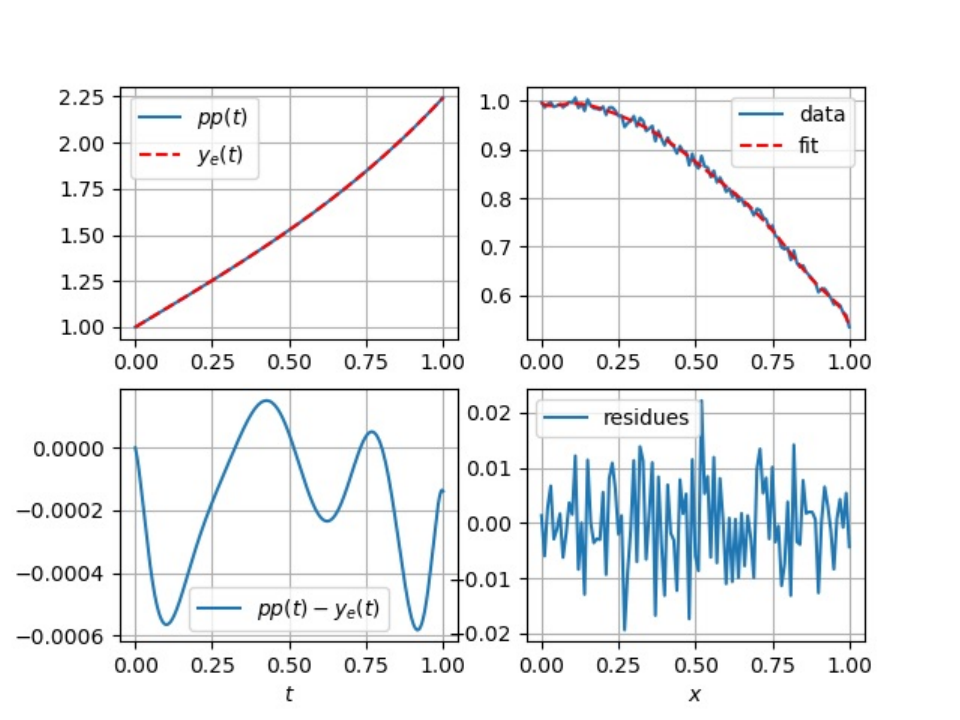}
\end{center}
\caption{ Solution of IDE (\ref{33}). $\sigma = 0.01$}
\label{Fig7}
\end{figure}

The graph in Fig. \ref{Fig7} shows the results of the numerical solution of task (\ref{34}). The right hand side of Eq.  (\ref{34}) has been computed on one hundred uniformly distributed points and contaminated by white noise with  $\sigma = 0.01$.  As can be seen from the graph, we have achieved excellent results.

\subsection{Linear integral equations of the second kind}
In this section, we will solve non-homogeneous Volterra   (\ref{V2}), and Fredholm  (\ref{Fr2})  integral equations of the second kind.

The Volterra integral equation  (\ref{V2}) has a unique solution for any continuous  function $\phi(x) $  \cite{Polyanin}. 
The Fredholm integral equation of the second kind  (\ref{Fr2})  has a unique solution if the parameter $\lambda$ is regular, meaning it is not equal to any eigenvalue of equation (\ref{Fr2})  \cite{Polyanin}.

\begin{equation}
\label{V2}
y(x) -  \int_a^x K(x, t)  y(t) \mathrm{d}t = \phi(x) 
 \end{equation}

\begin{equation}
\label{Fr2}
y(x) - \lambda \int_a^b K(x, t)  y(t) \mathrm{d}t = \phi(x)
 \end{equation}

 We assume kernels $K(x, t)$ and right hand side $\phi(x)$  of the equations under consideration have all properties necessary for the existence and uniqueness of the solution \cite{Polyanin}.
 
\subsubsection{Volterra integral equations of the second kind}

Let's solve the following Volterra integral equation of the second kind on interval $x \in [0, 5]$ using 100 points, presumably measured with $\sigma = 0.01$.

\begin{equation}
\label{35}
   y(x) + \int_0^x{(2 x - t) y(t) \mathrm{d}t} = -1 + \delta
 \end{equation} 
 
 The exact solution is $y(x) = (x^2 - 1) \mathrm{e}^{-x^2/2 }$. 

\begin{figure}[ht]
\begin{center}
\includegraphics[height=9.0 cm]{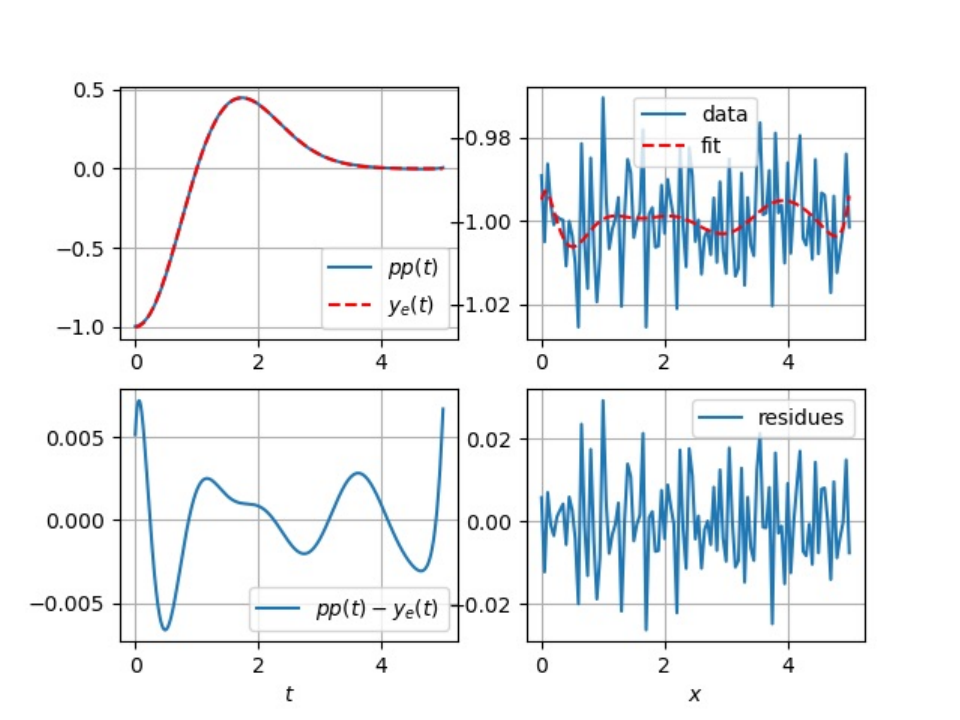}
\end{center}
\caption{ Solution of Volterra integral equation of the second kind (\ref{35}). $\sigma = 0.01$}
\label{Fig8}
\end{figure}

The piecewise-polynomial approximation to the exact solution of Eq. (\ref{35}) is shown in  Figure \ref{Fig8} . As can be seen from the graph, the sought function has been excellently approximated from discrete noisy data points.

\subsubsection{Fredholm integral equations of the second kind}

The  numerical solution of  Fredholm integral equations of the second kind (\ref{Fr2})  for regular values $\lambda$ is typically obtained by applying  the Gauss-Legendre quadrature rule for computing the integral. In a result, the solution is obtained only at the quadrature points,  and special care is needed for computing the solution at any point \cite{Press}.

 Let's solve the following Fredholm integral equation of the second kind, which has the exact solution $y(x) = \mathrm{e}^{-x}$.

\begin{equation}
\label{37}
   y(x) + \int_{0}^{b} \mathrm{e}^{-xt} y(t) \mathrm{d}t =    \frac{1 - \mathrm{e}^{-b(x + 1)}} {x+1}
 \end{equation} 
 
 For $ b = 100$ we have generated 100 noisy points with $\sigma = 10^{-4}.$

\begin{figure}[ht]
\begin{center}
\includegraphics[height=9.0 cm]{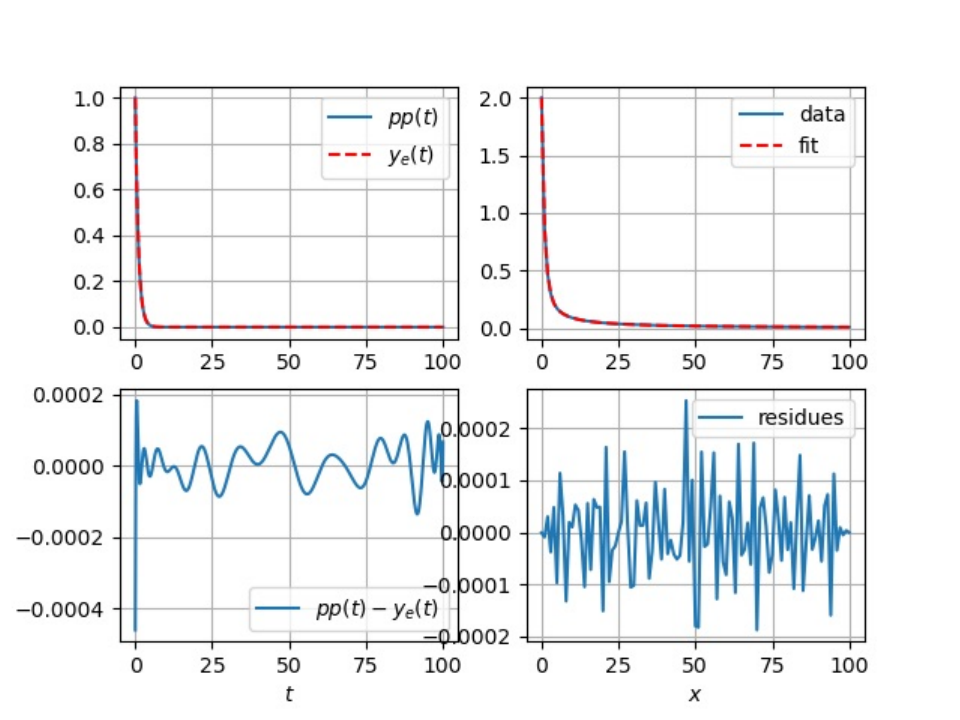}
\end{center}
\caption{ Solution of Fredholm integral equation  of the second kind (\ref{37}). $\sigma = 10^{-4}$}
\label{Fig9}
\end{figure}

The  solution of Fredholm integral equations of the second kind is stable with respect to changes in the right-hand side. This fact is confirmed by the graphs 
in (Figure \ref{Fig9}. As can be seen, the inaccuracy of the computed polynomial approximation and  the noise in the  right-hand side are of the same order.

\section{Solution of Ill-Posed Problems}

In general,  integral equations of the first kind are ill-posed, and regularization is usually required to obtain a stable solution.

\subsection{Volterra integral equations of the first kind}

Volterra integral equation of the first kind (\ref{V1}) 
\begin{equation}
\label{V1}
\int_a^x K(x, t)  y(t) \mathrm{d}t = \phi(x)
 \end{equation}

 has a unique solution under certain conditions \cite{Polyanin}, and is often solved by reducing it to Volterra equation of the second kind (\ref{V2}), or by using a method of quadratures  \cite{Polyanin}.  It's worth mentioning that differentiation of right-hand side function is required for reducing the equation to the second kind.

Here, we solve Volterra equations of the first kind for a noisy right hand side $\phi(x)$.  

\begin{equation}
\label{38}
   \int_0^x{\sin{(x - t)}y(t) \mathrm{d}t} = x \sin{x} + \delta
 \end{equation} 
 
 The exact solution of Eq. (\ref{38}),  $y(x) = \cos{x}$,  can be easily verified directly.
 
  \begin{figure}[ht]
\begin{center}
\includegraphics[height=9.0 cm]{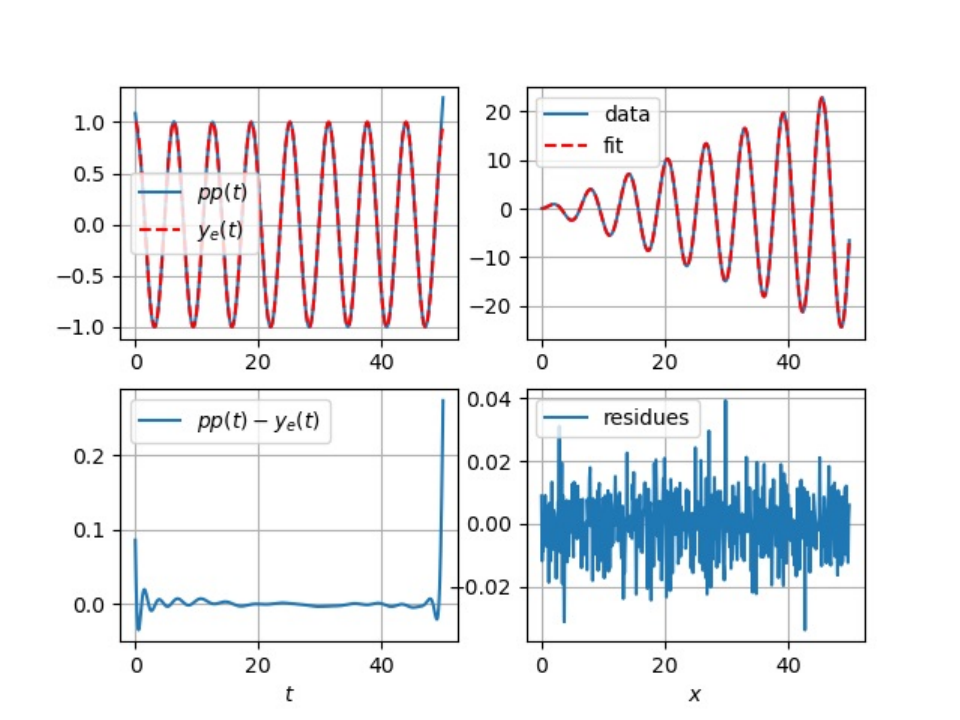}
\end{center}
\caption{ Solution of Volterra integral equation  of the first kind (\ref{38}). $\sigma = 10^{-2}$}
\label{Fig10}
\end{figure}

For numerical solution shown in Figure \ref{Fig10},  the right-hand side was computed for $m = 500$ points on the interval $[0, 50]$ with additive white noise with $\sigma = 10^{-2}$.  As can be seen from the graph, the obtained solution is slightly less accurate at the ends, but there is no error accumulation problem.

As next example consider the Abel integral equation.

\begin{equation}
\label{39}
   \int_0^x{\frac{y(t)}{\sqrt{x - t} }\mathrm{d}t} = \phi(x) + \delta
 \end{equation} 
 
 As can be seen, the kernel of Eq. (\ref{39}) has an integrable singular point at the upper integration limit,  and discretization of this equation requires special attention.
 For this example, we selected the exact solution $y(x) = \sin{x} \text{ on the interval } x \in [0, 3]$ and computed $\phi(x) $ numerically with the help of python's scipy.integrate.quad routine, adding  noise with $\sigma = 10^{-2}$.
 
 \begin{figure}[ht]
\begin{center}
\includegraphics[height=9.0 cm]{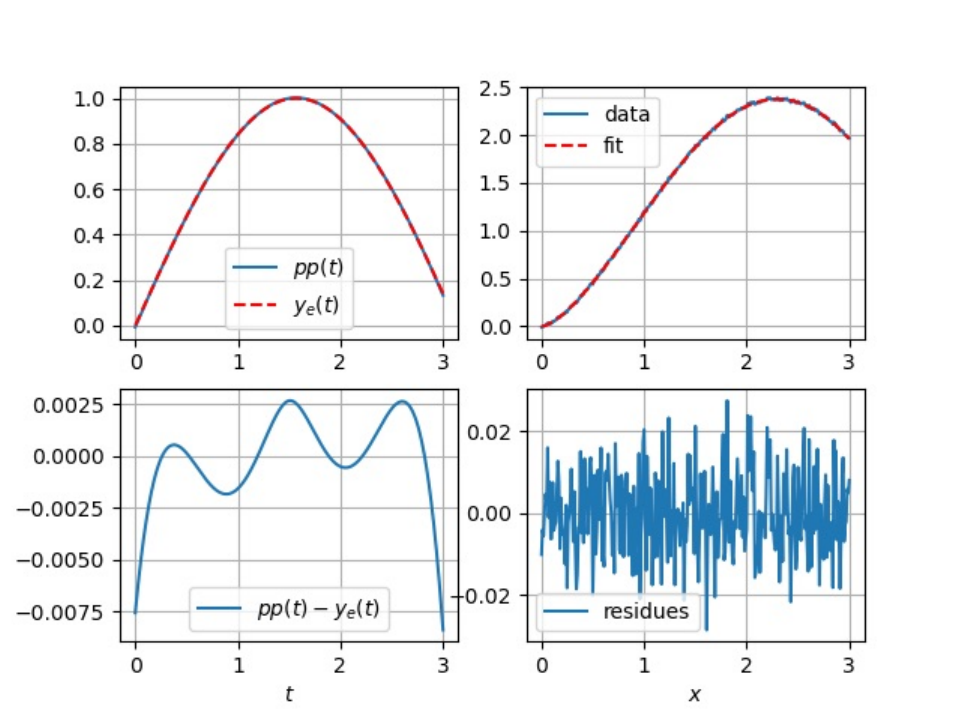}
\end{center}
\caption{ Solution of Abel integral equation. $\sigma = 10^{-2}$}
\label{Fig11}
\end{figure}

As shown in the graphs in Fig \ref{Fig11},  we have obtained excellent results for noisy right hand side with  $\sigma = 10^{-2}$.

\subsection{Fredholm integral equations of the first kind}

 \begin{equation}
\label{Fr1}
\int_a^b K(x, t)  y(t) \mathrm{d}t = \phi(x),  x \in [c, d]
 \end{equation}

It can be said that the Fredholm integral equations of the first kind are more challenging than Volterra integral equations.  

Let's  find a numerical solution of the following integral equation:

 \begin{equation}
\label{Fr1_1}
\int_0^{10}  y(t) \sin{xt}  \mathrm{d}t = \frac{\sin{10x} - 10x\cos{10x}}{x^2} + \delta, \\
  x \in [-1, 15]
 \end{equation}
 The exact solution of equation (\ref{Fr1_1}),  $y(x) = x$,  can be verified by direct integral evaluation.
 
 \begin{figure}[ht]
\begin{center}
\includegraphics[height=9.0 cm]{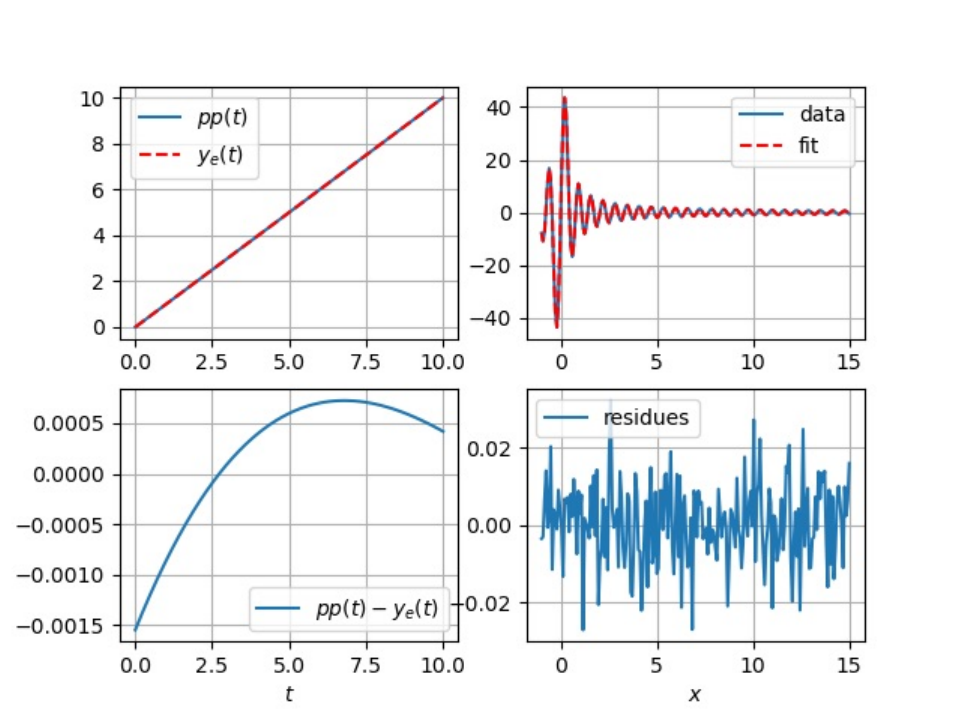}
\end{center}
\caption{ Regularized solution of equation  (\ref{Fr1_1}).  $\sigma = 10^{-2}$}
\label{Fig13}
\end{figure}

The results  shown in Figure \ref{Fig13} are likely as good as possible.
\\

Inverting  Laplace transforms given on the real axis is one of the most renowned challenging  problems  associated with Fredholm integral equations of the first kind with infinite upper limit of integration \cite{Varah}.

The damping property of the exponential kernel $\exp(-xt) $ significantly limits the information content available from  real-valued Laplace transforms.  The limitations of inverting real-valued Laplace transforms have been analyzed in    \cite{kr1}, and  can be summarized as follows:
\\

{A Laplace transform can be inverted from the real axis data if its theoretical inverse, when plotted on a semilogarithmic scale, is monotonic or has few extrema. The sharpness of resolvable peaks depends on the noise level in input data.}
\\

To facilitate the numerical inversion of the real Laplace transforms, data was synthesized using a combination of gamma distributions:

\begin{eqnarray}
\label{40}
\nonumber F(x)  = \sum_{i=1}^2{a_i}\left( \frac{\beta_i}{x + \beta_i} \right)^{\alpha_i} \\
\nonumber y(t) = \sum_{i=1}^2{a_i}\frac{t^{\alpha_i-1}\mathrm{e}^{-\beta_i t}\beta_i^{\alpha_i}}{\Gamma(\alpha_i)} \\
 \int_0^{\infty}{ {\mathrm{e}^{-x t} y(t)}\mathrm{d}t}  =  F(x) + \delta,
 \end{eqnarray} 
 
where $\Gamma(\alpha) $ is gamma function, and  parameters are $a = [1, 6] $, $\alpha = [10, 5]$, and $\beta = [50, 5]$.

 \begin{figure}[ht]
\begin{center}
\includegraphics[height=9.0 cm]{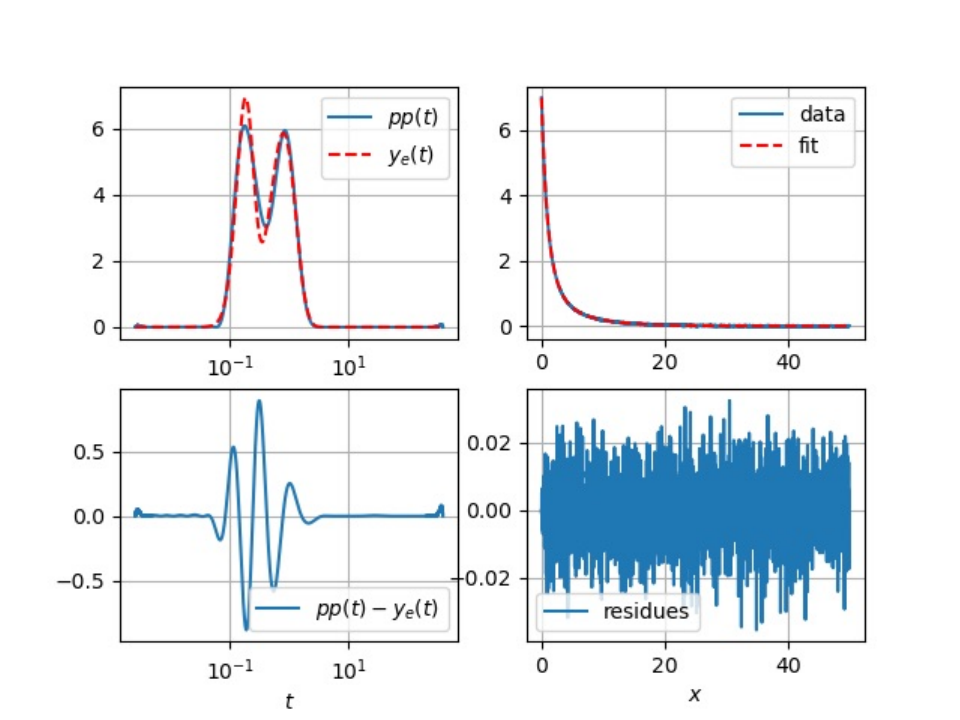}
\end{center}
\caption{ Inverse Laplace transform of a sum of gamma distributions (\ref{40}).  $\sigma = 10^{-2}$}
\label{Fig12}
\end{figure}

Graphs in Figure \ref{Fig12} show restoration of a two-peak function from its Laplace transform.  
As can be seen the peaks are restored quite accurately.

The problem of deconvolving nuclear magnetic resonance (NMR) data is closely related to inverse Laplace transformations and involves inverting the following transformation:
\begin{equation}
\label{41}
 \int_0^{\infty}{ {\mathrm{e}^{-x / t} y(t)}\mathrm{d}t}  =  F(x) + \delta.
 \end{equation} 
 
 For demonstration purposes, we use a sum of gamma distributions:
 \begin{equation}
 \nonumber y(t) = \sum_{i=1}^3{a_i}\frac{t^{\alpha_i-1}\mathrm{e}^{-\beta_i t}\beta_i^{\alpha_i}}{\Gamma(\alpha_i)}, 
  \end{equation} 
  with parameters $a = [1, 2, 6],  \alpha = [10,13, 15]$,  and  $ \beta=[100, 13, 1.5] $ .
  
  The right hand side $F(x)$ in Eq. (\ref{41}) was computed numerically, and contaminated data by white noise with $\sigma = 0.01$.

 \begin{figure}[ht]
\begin{center}
\includegraphics[height=9.0 cm]{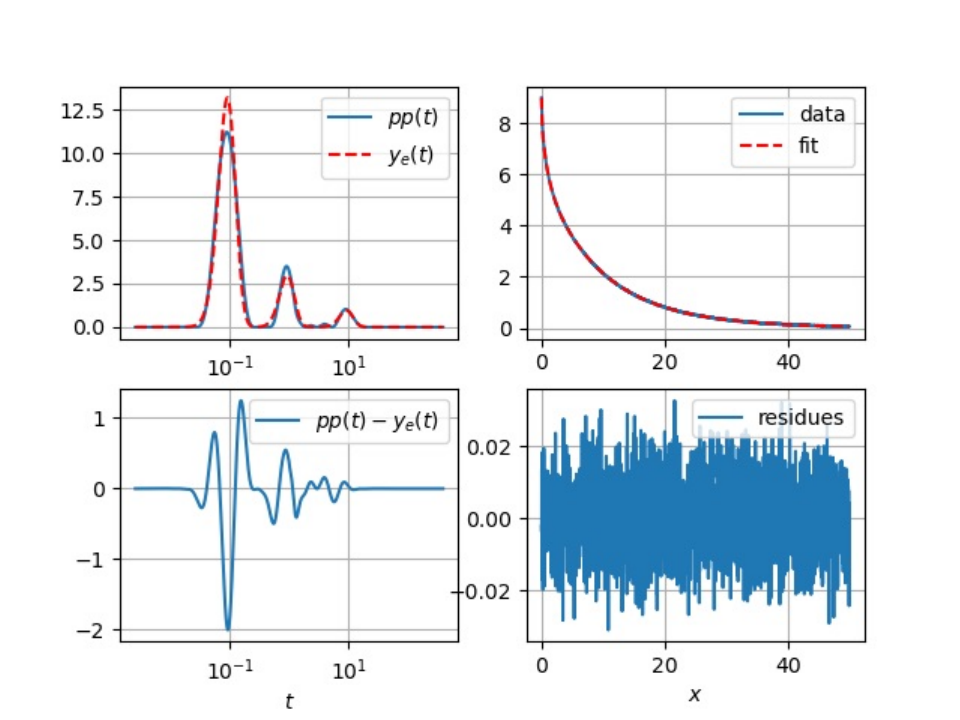}
\end{center}
\caption{NMR deconvolution.  $\sigma = 0.01$.}
\label{Fig14}
\end{figure}

The results of the NMR deconvolution are shown in Figure \ref{Fig14}. As can be seen from the graphs, all peaks were clearly resolved from quite noisy data.

Let's mention that positiveness of inverses shown in Figures \ref{Fig12} and \ref{Fig14} has been enforced.

\section{High-Degree Polynomial Interpolation \label{S51}}

High degree piecewise-polynomial approximation, considered in Section \ref{S2}, works well when the number of data points significantly exceeds the total number of polynomial coefficients, $ m \gg n$. If data is obtained from an experiment, it might not be possible to achieve this condition.

Cases where the function in question and the number of data points are such that $m \approx n$ are similar to high degree polynomial interpolation.
High degree polynomial interpolation is recommended to be avoided, as probably noted in every book on numerical analysis \cite{Kahaner, Press}. Because of ill-conditioning, high degree interpolating polynomials are not a good fit for sampled functions: "small changes in the data can give large differences in the oscillations between the points" \cite{Press}.
\\ 

As demonstrated in Section \ref{S4}, we determine the regularized polynomial coefficients by incorporating the stabilizing matrix   $ L $, which is computed for a specified number of arbitrarily distributed points. Subsequently, for  $ \lambda > 0 $  and $ k \gg m $, the second term in Equation (\ref{28}) should mitigate the large oscillations between data points.

 \begin{figure}[ht]
\begin{center}
\includegraphics[height=9.0 cm]{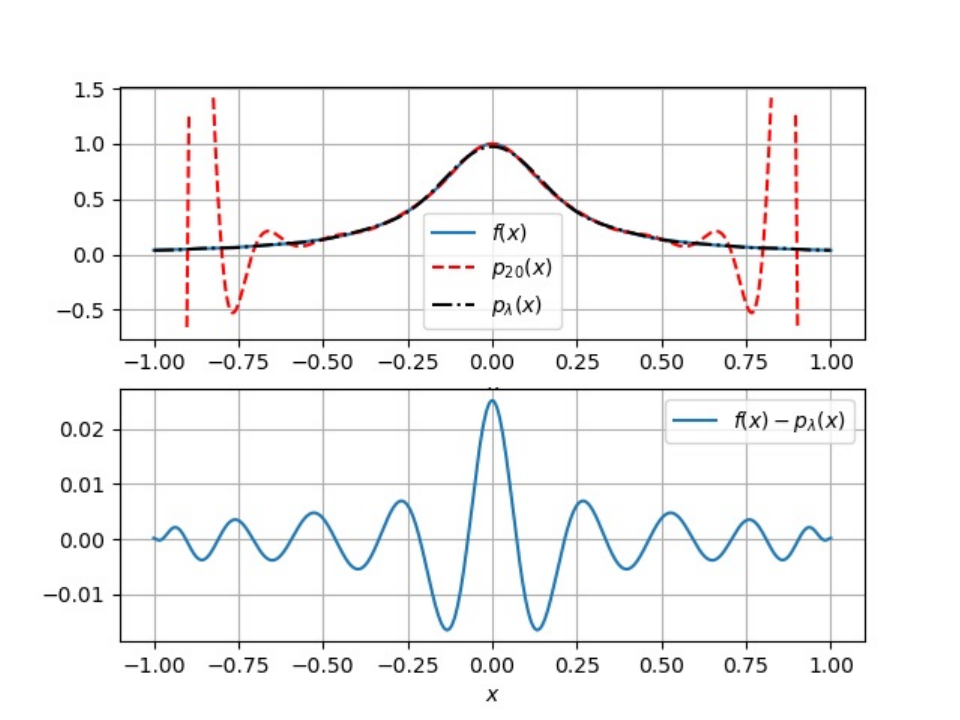}
\end{center}
\caption{Standard $p_{20}$ and stabilized $p_{\lambda}$ polynomial interpolations of Runge's function using $21$ uniformly distributed points. }
\label{Fig2}
\end{figure}

\begin{figure}[ht]
\begin{center}
\includegraphics[height=9.0 cm]{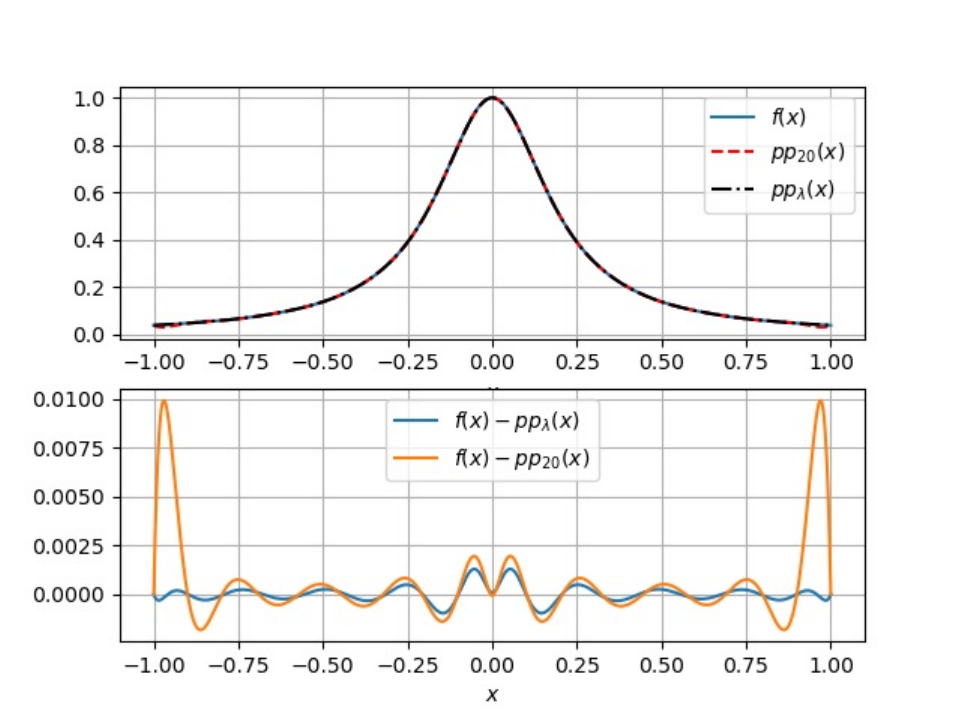}
\end{center}
\caption{ Piecewise-polynomial interpolations of Runge's function computed using 21 points. }
\label{Fig2a}
\end{figure}

The graph in Fig. \ref{Fig2} illustrates Runge's function $f(x) = \frac{1}{1 + 25x^2}$, its interpolating polynomial $p_{20}(x)$, and the stabilized interpolating polynomial $p_{\lambda}(x)$. The stabilizing matrix $L$ was computed using $k=1000$ uniformly distributed points. As evident from the graph, the addition of a stabilizing term indeed suppresses the large oscillations between the points.

Furthermore, Fig. \ref{Fig2a} displays piecewise-polynomial interpolations of Runge's function computed using 21 points, with $pp_{\lambda}(x)$ including the stabilizing term, and without, represented by $pp_{20}(x)$. It's observable that the large oscillations have been eliminated in both cases. However, the piecewise-polynomial computed with the stabilizing term is more accurate, particularly at the ends.

The introduction of the stabilizing matrix $ L $ plays a crucial role in regularization, here we applyed it for stabilizing a well-posed yet ill-conditioned task, thereby improving the accuracy of the polynomial interpolation.

Let us  note that the ability to solve a problem both with and without the stabilizing term could  be utilized to test the stability of the specific problem.

\section{Conclusion}

In this paper, we have explored solutions to Eq. (\ref{4}) across various special cases. By representing the approximate solution of Eq. (\ref{4}) in the form of high-degree piecewise-polynomials, we have demonstrated that:

\begin{itemize}
    \item Problems involving differentiation and integration can be uniformly addressed by transforming them into a minimization problem.
    \item The unified numerical scheme is readily adaptable for solving ill-posed problems.
    \item In certain cases, well-posed problems may benefit from the addition of a stabilizing term, similar to regularization.
\end{itemize}

\section{Motivation and afterword \label{A}}

\epigraph{" you have to work hard to get your thinking clean to make it simple. "} {Steve Jobs: There's Sanity Returning}

The suggested numerical approach seems self-evident, though it has been overlooked until now. In this section, the author's motivation and the steps that led to the proposed numerical scheme are briefly discussed.\\

The regularized solution of linear integral equations of the first kind holds a central position in regularization theory \cite{Tikhonov}. However, it has been recognized that conventional regularization methods are insufficient for inverting Laplace transforms given on the real axis \cite{Varah}.

An integral form of the regularizing operator for inverting Laplace transforms computable at any point on the real axis has been analytically derived by the author \cite{kr1}. This regularizing operator differs from conventional regularizing operators in that it has two tuning parameters along with the regularization parameter. The accuracy of the computed inverse significantly depends on the values of the tuning parameters.

Optimal values of all parameters were computed with the help of a proposed ad hoc criterion for finding the regularization parameter \cite{kr2}.

It can be argued that this method allows for extracting more information from given data than any other real-valued method. One can test this claim using software available on GitHub \cite{github}.

Subsequently, the effort to develop parametrized regularizing operators for any Fredholm integral equation of the first kind \cite{arxiv}, prompted the representation of the sought function $y(x)$ as a piecewise-polynomial. This approach resulted in a matrix equation whose regularization is somewhat nonstandard.

Following the regularization of this equation (as discussed in Section \ref{S4}), the piecewise-polynomial approximation of the solution was employed to solve Fredholm integral equations of the second kind. 

The solution of Fredholm integral equations of the second kind has shown that parametrized inverse operators work well for well-posed problems, and their applications are not limited to ill-posed problems.  Consequently, integro-differential and ordinary differential equations have also been included in the consideration.\\

It may also be noteworthy to mention that the author's research has never been purely academic but has always been pursued during free time, especially throughout the years of retirement. This freedom to explore seemingly unconventional directions has proven invaluable.
\\



\end{document}